\documentclass{amsart}
\newtheorem{theorem}{Theorem}

\theoremstyle{definition}

\theoremstyle{remark}

\numberwithin{equation}{section}
\theoremstyle{corollary}

\theoremstyle{proposition}
\newtheorem{proposition}[theorem]{Proposition}
\newfont{\EUL}{eufm10 scaled 1000}

\newcommand\Lt{\mbox{\EUL t}}

\renewcommand\k{\mbox{\EUL k}}

\newcommand\m{\mbox{\EUL m}}

\newcommand\Ad{{\rm Ad}}

\begin{document}
%
%
\title{A note on the moment map on compact K\"ahler manifolds}
\author{Anna Gori and Fabio Podest\`a}
\address{Dipartimento di Matematica \lq U.Dini\rq\\ Viale Morgagni 67/A\\
50100 Firenze\\Italy} \email{gori@math.unifi.it}
\address{Dipartimento di Matematica e Appl.\ per l'Architettura\\
Piazza Ghiberti 27\\50100 Florence\\Italy} \email{podesta@math.unifi.it} \subjclass{}
\begin{abstract}
We consider compact K\"ahler manifolds acted on by a connected
compact Lie group $K$ of isometries in a Hamiltonian fashion. We prove
that the squared moment map $||\mu||^2$ is constant if and only if
the manifold is biholomorphically and $K$-equivariantly isometric
to a product of a flag manifold and a compact K\"ahler manifold
which is acted on trivially by $K$.\\
The authors do not know whether the compactness of $M$ is
 essential in the main theorem;
more generally it would
be interesting to have a similar result for (compact) symplectic manifolds.
\end{abstract}
\maketitle
\section{Introduction}
We shall consider compact K\"ahler manifolds $M$ acted on by a
compact connected Lie group $K$ of isometries; such isometries are
automatically holomorphic transformations of $M$. We shall also
suppose that the $K$-action on $M$ is Hamiltonian, i.e. there
exists a moment map $\mu: M \to \k^*$, where $\k$ is the Lie
algebra of $K$; such a moment map exists if and only if the Lie
group $K$ acts trivially on the Albanese torus $Alb(M)$ (see
\cite{HW},\cite{GS}). Throughout the following we will denote by
$g$, $J$ and $\omega$ the Riemannian metric, the complex structure
and the corresponding K\"ahler form on $M$ respectively; moreover Lie groups
and their Lie algebras will be indicated with capital and gothic letters respectively. \\
If we fix an $\Ad(K)$-invariant scalar product $q:= \langle,\rangle$ on $\k$
and we identify $\k^*$ with $\k$ by means of $q$, we can think of $\mu$ as a
$\k$-valued map ;the function $f\in C^\infty(M)$ defined as $f:= ||\mu||^2$ has been extensively
used in \cite{Kir1} to obtain strong information on the topology of the manifold. \\
Our main result is the following
\begin{theorem}\label{splitting} Suppose
 $M$ is a compact K\"ahler $K$-Hamiltonian manifold, where
$K$ is a compact connected Lie group of isometries.
If $\mu$ denotes the moment map, then the squared moment map $f=\|\mu\|^2$
is constant if and only if the manifold $M$ is biholomorphically and $K$-equivariantly
isometric to the product of a flag manifold and a compact K\"ahler manifold
which is acted on trivially by $K$.
\end{theorem}
In order to prove the above theorem, we need the following result, which has been proved in \cite{GP}.
\begin{proposition}\label{max} Let $M$ be a compact K\"ahler manifold which is
acted on by a compact connected Lie group $K$ of isometries in a
Hamiltonian fashion with moment map $\mu$. If a point $x\in M$
realizes the maximum of $||\mu||^2$, then the orbit $K\cdot x$ is
complex.\end{proposition}
The authors do not know whether the compactness of $M$ is essential in the main theorem;
more generally it would
be interesting to have a similar result for (compact) symplectic manifolds.\\
\section{Proof of the main result}
For later use, we reproduce here the proof of Proposition \ref{max}.
\begin{proof} [Proof of Proposition \ref{max}] We will follow the notations as in \cite{Kir1}.
Let $\beta = \mu(x)$, which
we can suppose to lie in the closure of a Weyl chamber $\Lt_+$, where $\Lt$ denotes
the Lie algebra of a fixed maximal torus in $K$. If $\mu_\beta := \langle\mu,\beta\rangle$
is the height function relative to $\beta$ and if $Z_\beta$ denotes the union of the
connected components of the critical point set of $\mu_\beta$ on which $\mu_\beta$ takes the
value $||\beta||^2$, then $x$ belongs to the critical set $C_\beta =
K\cdot(Z_\beta\cap \mu^{-1}(\beta))$. We now claim that $Z_\beta =
\mu^{-1}(\beta)$. Indeed, if $p\in Z_\beta$, then $\mu_\beta(p) =
||\beta||^2$ and
$$||\beta||^2 \leq \langle\mu(p),\beta\rangle \leq ||\mu(p)||\cdot
||\beta|| \leq ||\beta||^2,$$ and therefore $\mu(p) = \beta$, i.e.
$p\in \mu^{-1}(\beta)$. Viceversa, if $p\in \mu^{-1}(\beta)$, then
$||\mu(p)||^2$ is the maximum value of $f:= ||\mu||^2$ and
therefore $\hat{\beta}_p = 0$, where $\hat\beta$ denotes the Killing field on $M$
induced by the element $\beta\in \k$; moreover $\mu_\beta(p) =
||\beta||^2$ and therefore $p\in Z_\beta$. This implies that
$\mu^{-1}(\beta)$ is a complex
submanifold and that $C_\beta = K\cdot \mu^{-1}(\beta)$. \\
If $S_\beta$ denotes the Morse stratum of $C_\beta$, we claim that $S_\beta = C_\beta$.
Indeed, if $\gamma_t(q)$
denotes the flow of $-grad(f)$ through a point $q$ belonging to
the stratum $S_\beta$, then $\gamma_t(q)$ has a limit point in the
critical subset $C_\beta$; since $f(\gamma_t(q))$ is
non-increasing for $t\geq 0$ and $f(C_\beta)$ is the maximum value
of $f$, we see that $f(\gamma_t(q)) = ||\beta||^2$ for all $t\geq
0$, that is $S_\beta\subseteq C_\beta$ and therefore
$S_\beta = C_\beta$.\\
This implies that $C_\beta = S_\beta$ is a smooth complex
submanifold of $M$ and for every $y\in
\mu^{-1}(\beta)$, we have
$$T_yS_\beta = T_y(K\cdot y) + T_y(\mu^{-1}(\beta)).$$
Now, if $v\in T_y(\mu^{-1}(\beta))$, then $v = Jw$ for some $w\in
T_y(\mu^{-1}(\beta))$ and for every $X\in \k$ we have
$$0 = \langle d\mu_y(w),X\rangle = \omega_y(w,{\hat X}_y) =
\omega_y(Jv,{\hat X}_y) = g_y(v,{\hat X}_y),$$ meaning that
$T_y(\mu^{-1}(\beta))$ is $g$-orthogonal to $T_y(K\cdot y)$. Since
both $S_\beta$ and $\mu^{-1}(\beta)$ are complex, this implies
that $K\cdot y$ is a complex orbit.\end{proof}
We now give the proof of Theorem \ref{splitting}

\begin{proof}[Proof of Theorem \ref{splitting}]
Assume $f$ to be constant, i.e. the manifold $M$ is
mapped, by $\mu$, into a sphere. We fix a point $x_o\in M$;
we also recall that $\mu(M)\cap
\mathfrak{t}^*_+$ is convex \cite{Kir2}, hence the manifold is mapped to a
single coadjoint orbit $\mathcal{O}$ $=K/K_{\beta}$, where $\beta = \mu(x_o)$.
We then have that all the points of $M$ are critical for
$f$ and $M=K\cdot \mu^{-1}(\beta)$ by the arguments used in the proof
of Proposition ~\ref{max}.\\
Note that each $K$-orbit
intersects $\mu^{-1}(\beta)$ in a single point. Indeed, if
there are two points $x$ and $z=k\cdot x$ for $k\in K$ which
lie in $\mu^{-1}(\beta)$, then, by
the $K$-equivariance of $\mu,$ we have $k\in K_{\beta}$ which
is equal to $K_x$, since, for Proposition \ref{max}, the $K\cdot
x$ orbit is complex; hence $k\cdot x = x$.\\
We can also prove that, for each $x\in M,$
the tangent space $T_x K\cdot x$ is orthogonal to $T_x
\mu^{-1}(\mu(x))$; indeed if $ X\in T_x K\cdot x$ and $Y\in
T_x \mu^{-1}(\mu(x))$, then, using the fundamental property of the
moment map $\mu$, we argue that $0=\omega(Y, X)=g(Y,J X)$,
where $JX \in T_x K\cdot x$ since the orbit is complex.\\
From this it follows that the map
$$\varphi: K/K_\beta\times \mu^{-1}(\beta)\to M,\qquad \varphi(gK_\beta,x) = g\cdot x,$$
where we identify $K/K_\beta$ with the orbit $K\cdot x_o$, is a well defined $K$-equivariant
diffeomorphism.\\
We also observe that $\mu^{-1}(\mu(x))$ is connected for all $x\in M$;
moreover all the $K$-orbits
are principal since their stabilizers are all equal to
$K_{\beta}$, hence we have
that $K_x$ acts trivially on $T_x \mu^{-1}(\mu(x))=(T_x K\cdot
x)^\perp$ for all $x\in M$. \\
We now denote by $\mathcal F$ the foliation given by the $K$-orbits and by
${\mathcal F}^\perp$ the orthogonal foliation, so that ${\mathcal F}^\perp_y =
T_y(\mu^{-1}(\mu(y)))$. We now claim that both $\mathcal F$ and ${\mathcal F}^\perp$
are totally geodesic; this will then imply that they are both parallel with respect to the Levi
Civita connection and our result will follow. \\
We first observe that ${\mathcal F}^\perp$ is totally geodesic. Indeed,
at each point $y\in M$, the stabilizer $K_y$ is the centralizer of some torus in $K$ and therefore
its isotropy representation on the tangent space $T_y(K\cdot y)$ has no fixed vector; on the other
hand, $y$ is a principal point, so that the isotropy representation of $K_y$ leaves
${\mathcal F}^\perp_y$ pointwisely fixed. This shows that $\mu^{-1}(\mu(y))$ is a connected
component of the fixed point set of $K_y$ and therefore it is totally geodesic.\\
We now claim that $\mathcal F$ is totally geodesic. We fix again a point $y\in M$
and a normal vector $\xi \in T_y(K\cdot y)^\perp$; the shape operator
$A_\xi$ of the orbit $K\cdot y$ relative to the normal vector $\xi$ is
a self-adjoint operator on $T_y(K\cdot y)$, which is $K_y$-invariant, since $K_y$ leaves $\xi$ fixed.
We now decompose the Lie algebra $\k$ of $K$ as
$\k = {\k}_\beta \oplus \m,$ where ${\k}_\beta$ is the Lie algebra of $K_\beta$ and
$\m$ is an $\operatorname{Ad}(K_\beta)$-invariant subspace which can be
identified with the tangent space $T_y(K\cdot y)$. We also split
$\m = \bigoplus_{i=1}^l {\m}_i$ as a sum of
$\operatorname{Ad}(K_\beta)$-irreducible submodules; it is known that
the summands ${\m}_i$ are non-trivial and mutually inequivalent as
$\operatorname{Ad}(K_\beta)$-modules (see e.g.~\cite{S}). This means that $A_\xi$ preserves each
${\m}_i$ and therefore, by Schur's Lemma, its restriction on each
${\m}_i$ is a multiple of the identity. The complex structure $J$, viewed as a
$\operatorname{Ad}(K_\beta)$-invariant operator on $\m$ also preserves
each submodule and therefore it commutes with $A_\xi$. \\
On the other hand, it is known (see e.g.~\cite{KN},~p. 175) that the shape operator of a complex
submanifold in a K\"ahler manifold anti-commutes with the complex
structure; from this we conclude that $A_\xi = 0$ for all normal vectors
$\xi$. This means that $\mathcal F$ is totally geodesic.\\
It is now easily seen that two orthogonal, integrable and totally geodesic
foliations are parallel and this concludes the proof.
\end{proof}

\end{document}